\documentclass[12pt]{amsart}
\usepackage{amsxtra,latexsym,amssymb}
\usepackage{epsfig,rotate,amsthm}
\usepackage{hyperref}

\def\qed{{\hfill $\Box$}}

\def\Z{{\mathbb Z}}
\def\C{{\mathbb C}}

\def\U{{U_{r,s}^{+}(B_{2})}}
\def\V{{\check U}^{\geq 0}_{r,s}(B_{2})}
\theoremstyle{theorem}
\newtheorem{thm}{Theorem}[section]
\newtheorem{cor}{Corollary}[section]
\newtheorem{prop}{Proposition}[section]

\newtheorem{lem}{Lemma}[section]
\theoremstyle{definition}
\newtheorem{defn}{Definition}[section]
\theoremstyle{remark}

\begin{document}
\title[Derivations of $\U$]{Derivations of the two-parameter quantized enveloping algebra $\U$}
\author[X. Tang]{Xin Tang}
\address{Department of Mathematics \& Computer Science\\
Fayetteville State University\\
1200 Murchison Road, Fayetteville, NC 28301}
\email{xtang@uncfsu.edu} 
\thanks{This research project is partially supported by a faculty research mini-grant from HBCU Master's Degree STEM Program at Fayetteville State University}
\date{\today}
\begin{abstract}
Let $r,s$ be two parameters chosen from $\C^{\ast}$ such that $r^{m}s^{n}=1$ implies $m=n=0$. We compute the derivations of the two-parameter quantized enveloping algebra $\U$ and calculate its 
first degree Hochschild cohomology group. We further determine the group of algebra
automorphisms for the two-parameter Hopf algebra $\V$. As a result, we determine 
the group of Hopf algebra automorphisms for $\V$. 
\end{abstract}
\maketitle 
\section*{Introduction}

Let $\mathfrak g$ be a finite dimensional complex simple Lie algebra. The two-parameter quantized enveloping algebras ( or quantum groups) $U_{r,s}({\mathfrak g})$ have been studied by Benkart and Witherspoon in the references \cite{BW1, BW2} for Lie algebras of type $A$. The two-parameter quantized enveloping algebras $U_{r,s}(\mathfrak g)$ have been further studied  for Lie algebras $\mathfrak g$  of type $B, C, D$ in \cite{BGH}. Overall, the two-parameter quantized enveloping algebras $U_{r,s}(\mathfrak g)$ are close analogues of their one-parameter peers, and share a similar structural and representation theory as the one-parameter quantized enveloping algebras $U_{q}(\mathfrak g)$. For instance, the positive part of the two-parameter quantized enveloping algebra $U_{r,s}(\mathfrak g^{+})$ is also proved to be isomorphic to certain two-parameter Ringel-Hall algebra in \cite{T1}. 

Nonetheless, there are some differences between the one-parameter quantized enveloping algebras and two-parameter quantized enveloping algebras. In particular, the center of two-parameter quantized enveloping algebras $U_{r,s}(\mathfrak g)$ already posed a different 
picture \cite{BKL}.  Besides, it seems that the two-parameter quantized enveloping algebras $U_{r,s}(\mathfrak g)$ are more rigid in that they possess less symmetry. On the one hand, 
these differences make it more interesting to further investigate these algebras. On the other hand, these differences also make it plausible  to more effectively study the structures of 
these algebras.

Recently, there have been some interests in the study of the subalgebras  of two-parameter quantized enveloping algebras. For example, both the derivations of the subalgebra $U_{r,s}^{+}(sl_{3})$ and the automorphisms of augmented Hopf algebra $\check{U}_{r,s}^{\geq 0}(sl_{3})$ have been determined in \cite{T2}. In this paper, we will study the derivations of the algebra 
$\U$ and the automorphisms of the Hopf algebra $\V$. First of all, we will completely determine 
the derivations and calculate the first degree Hochschild cohomology group for the algebra $\U$. Second of all, we will determine both the algebra automorphism group and Hopf algebra automorphism group of the Hopf algebra $\V$. To calculate the derivations, one needs to embed the algebra $\U$ into a quantum torus, where the derivations are known \cite{OP}. Via this embedding, we shall be able to pull the information on derivations back to the algebra $\U$. As a matter of fact, applying a result in \cite{OP}, we will be able determine all the derivations of the algebra $\U$ up to its inner derivations. As a result, we show that the first degree Hochschild cohomology group $HH^{1}(\U)$ of $\U$ is indeed a $2-$dimensional vector space over the base field $\C$. In order to determine the automorphisms of $\V$, we will follow the lines in \cite{F}; and we obtain similar results as those in \cite{F,T2}.

The paper is organized as follows. In Section 1, we recall some basic definitions and properties on the two-parameter quantized enveloping algebras $\U$, and establish some necessary commuting identities; then we determine the derivations and calculate the first degree Hochschild cohomology group. In Section 2, we first determine the algebra automorphism group of $\V$; and then we determine the Hopf algebra automorphism group of $\V$. 

\section{Derivations and the first degree Hochschild Cohomology group of $\U$}

In this section, we compute all the derivations of the two-parameter quantized enveloping algebra $\U$. As a matter of fact, we are able to show that every derivation of $\U$ can be uniquely decomposed as the sum of an inner derivation plus a linear combination of certain specifically defined derivations. As a result, we prove that the first degree Hochschild cohomology group of $\U$ is a two-dimensional vector space over the center of $\U$, which can be proved to be the base field $\C$. The computation of derivations will be carried out via an embedding of the algebra $\U$ into a quantum torus, whose derivations had been described in \cite{OP}; and this embedding allows to pull information on derivations of the quantum torus back to the algebra $\U$. We should mention that this method has also been successfully used to compute the derivations of some quantum algebras such as $U_{q}({\mathfrak sl_{4}^{+}})$ in \cite{LL} and $U_{r,s}^{+}(sl_{3})$ in \cite{T2}. In addition, the derivations of $U_{q}^{+}(B_{2})$ were determined in \cite{BL} via a similar approach.

\subsection{Some basic properties of the algebra $\U$}

Recall that two-parameter quantum groups $U_{r,s}(B_{n})$ associated to the complex simple Lie algebras of type $B_{n}, n\geq 2$ have been studied in \cite{BGH}. For our convenience, we will recall the construction for one of the subalgebras, i.e., the algebra $\U$ in the case $B_{2}$ here.  From now on, we will always assume that the parameters $r, s$ are chosen from $\C^{\ast}$ such that $r^{m}s^{n}=1$ implies $m=n=0$. 

First of all, we recall the following definition from the reference \cite{BGH}:
\begin{defn}
[See \cite{BGH}] The two-parameter quantized enveloping algebra $\U$ is defined to be the $\C-$algebra generated by the generators $e_{1}, e_{2}$ subject to the following relations:
\begin{eqnarray*}
e_{1}^{2}e_{2}-(r^{2}+s^{2})e_{1}e_{2}e_{1}+r^{2}s^{2}e_{2}e_{1}^{2}=0,\\
e_{1}e_{2}^{3}-(r^{2}+rs+s^{2})e_{2}e_{1}e_{2}^{2}
+rs(r^{2}+rs+s^{2})e_{2}^{2}e_{1}e_{2}-r^{3}s^{3}e_{2}^{3}e_{1}=0.
\end{eqnarray*}
\end{defn}

In the rest of this subsection, we will establish some basic properties of the algebra $\U$. In particular, we will show that the two-parameter quantized enveloping algebra $\U$ can be presented as an iterated skew polynomial ring. As a result, we will be able construct a PBW-basis for $\U$, and prove that the center of $\U$ is reduced to $\C$. 

Now we fix some notation by setting the following new variables $X_{1}, X_{2}, X_{3}$ and $X_{4}$: 
\begin{eqnarray*}
X_{1}=e_{1},\quad X_{2}=e_{3}=e_{1}e_{2}-r^{2}e_{2}e_{1},\\
X_{3}=e_{2}e_{3}-s^{-2}e_{3}e_{2},\quad
X_{4}=e_{2}.
\end{eqnarray*}

Concerning the relations between these new variables, we shall have the following lemma:
\begin{lem} The following identities hold:
\begin{enumerate}
\item $X_{1}X_{2}=s^{2}X_{2}X_{1}$;\\
\item $X_{1}X_{3}=r^{2}s^{2}X_{3}X_{1}$;\\
\item $X_{2}X_{3}=rsX_{3}X_{2}$;\\
\item $X_{1}X_{4}=r^{2}X_{4}X_{1}+X_{2}$;\\
\item $X_{2}X_{4}=s^{2}X_{4}X_{2}-s^{2}X_{3}$;\\
\item $X_{4}X_{3}=r^{-1}s^{-1}X_{3}X_{4}$.
\end{enumerate}
\end{lem}
{\bf Proof:} These identities can be verified via straightforward computation and we will skip the details here.
\qed

In addition, let us define some algebra automorphisms $\tau_{2},\tau_{3}$, and $\tau_{4}$, and some derivations $\delta_{2},\delta_{3}$, and $\delta_{4}$ as follows: 
\begin{eqnarray*}
\tau_{2}(X_{1})=s^{-2}X_{1},\quad \delta_{2}(X_{2})=0,\\
\tau_{3}(X_{1})=r^{-2}s^{-2}X_{1},\,\tau_{3}(X_{2})=r^{-1}s^{-1}X_{2},\\
\delta_{3}(X_{1})=0,\,\delta_{3}(X_{2})=0,\\
\tau_{4}(X_{1})=r^{-2}X_{1},\,\tau_{4}(X_{2})=S^{-2}X_{2},\,\tau_{4}(X_{3})=r^{-1}s^{-1}X_{3},\\
\delta_{4}(X_{1})=-r^{-1}X_{2},\, \delta_{4}(X_{2})=X_{3},\, \delta_{4}(X_{3})=0.
\end{eqnarray*}
\qed

Based on the previous lemma, it is easy to see that we have the following result
\begin{thm}
The algebra $\U$ can be presented as an iterated skew polynomial ring. In
particular, we have the following result
\[
\U \cong
\C[X_{1}][X_{2},\tau_{2},\delta_{2}][X_{3},\tau_{3},\delta_{3}][X_{4}, \tau_{4}, \delta_{4}].
\] 
\end{thm}
\qed

Based on the previous theorem, we have an obvious corollary as follows:
\begin{cor}
The set 
\[
\{X_{1}^{a}X_{2}^{b}X_{3}^{c}X_{4}^{d}|a, b, c, d \in \Z_{\geq 0}\}
\]
forms a PBW-basis of the algebra $\U$. In particular, $\U$ has a $GK-$dimension 
of $4$.
\end{cor}
\qed

Associated to this iterated skew polynomial ring presentation of $\U$, one can define a filtration of $\U$ such that the corresponding graded algebra $gr\U$ is a quantum space generated by the variables $\overline{X_{1}}, \overline{X_{2}},\overline{X_{3}}$ and $\overline{X_{4}}$ subject to the following relations:
\begin{eqnarray*}
\overline{X_{1}}\,\overline{X_{2}}=s^{2}\overline{X_{2}}\,\overline{X_{1}};\\
\overline{X_{1}}\,\overline{X_{3}}=r^{2}s^{2}\overline{X_{3}}\,\overline{X_{1}};\\
\overline{X_{1}}\,\overline{X_{4}}=r^{2}\overline{X_{4}}\, \overline{X_{1}};\\
\overline{X_{2}}\,\overline{X_{3}}=rs\overline{X_{3}}\, \overline{X_{2}};\\
\overline{X_{2}}\, \overline{X_{4}}=s^{2}\overline{X_{4}}\, \overline{X_{2}};\\
\overline{X_{3}}\, \overline{X_{4}}=rs\overline{X_{4}}\,\overline{X_{3}}.\\
\end{eqnarray*}

Now we have the following description of the center of the algebra $\U$:
\begin{cor}
The center of $\U$ is reduced to the base field $\C$.
\end{cor}
{\bf Proof:} Let $u \in \U$ be an element in the center of $\U$. Then $u$ is a linear combination of the monomials $X_{1}^{a}X_{2}^{b}X_{3}^{c}X_{4}^{d}$. Then $\overline{u}$ is in the center of $gr\U$. Let $\overline{X_{1}^{a}}\,\overline{X_{2}^{b}}\,\overline{X_{3}^{c}}\,\overline{X_{4}^{d}}$ be the image of one of these monomials in the graded algebra, then this monomial  $\overline{X_{1}^{a}}\,\overline{X_{2}^{b}}\,\overline{X_{3}^{c}}\,\overline{X_{4}^{d}}$ commutes with the generators $\overline{X_{1}},\overline{ X_{4}}$. Therefore, we have the following
\begin{eqnarray*}
s^{2b}r^{2c}s^{2c}r^{2d}=1;\\
r^{-2a}s^{-2b}(rs)^{-c}=1.
\end{eqnarray*}
Therefore, we have the following
\begin{eqnarray*}
2b+2c=0,\quad 2c+2d=0;\\
2a+c=0,\quad
2b+c=0.
\end{eqnarray*}
Solving this system, we get $a=b=c=d=0$. Therefore, we have $u\in \C$. 
\qed

\subsection{The embedding of $\U$ into a quantum torus}
In this subsection, we construct an embedding of the algebra $\U$ into a quantum torus. This embedding shall enable us to extend the derivations of $\U$ to derivations of the quantum torus, and later on pull information backward. The point is that that the algebra $\U$ has a Goldie quotient ring, which we shall denote by $Q(\U)$. Within the Goldie quotient ring $Q(\U)$ of $\U$, let us define the following new variables
\begin{eqnarray*}
T_{1}=X_{1},\quad
T_{2}=X_{2},\\
T_{3}=X_{3},\quad
T_{4}=X_{2}^{-1}Z^{\prime}X_{1}^{-1},
\end{eqnarray*}
where 
\[
Z^{\prime}=(X_{1}(X_{3}+(s^{-2}-r^{-1}s^{-1})X_{2}X_{4})-s^{4}(X_{3}+(s^{-2}-r^{-1}s^{-1})X_{2}X_{4})X_{1}).
\]

Let us set a new variable $W=X_{3}+(s^{-2}-r^{-1}s^{-1})X_{2}X_{4}$, then we have the following lemma:
\begin{lem} The following identities hold:
\begin{enumerate}
\item $X_{1}W=r^{2}s^{2}WX_{1}+(1-r^{-1}s)X_{2}^{2}$;\\
\item $X_{2}W=s^{2}WX_{2}$;\\
\item $X_{3}W=WX_{3}$;\\
\item $X_{4}W=s^{-2}WX_{4}$;\\
\item $X_{1}Z^{\prime}=r^{2}s^{2}Z^{\prime}X_{1}$;\\
\item $X_{2}Z^{\prime}=Z^{\prime}X_{2}$;\\
\item $X_{3}Z^{\prime}=r^{-2}s^{-2}Z^{\prime}X_{3}$;\\
\item $X_{4}Z^{\prime}=r^{-2}s^{-2}Z^{\prime}X_{4}.$
\end{enumerate}
\end{lem}
{\bf Proof:} Once again, we can verify these identities by brutal force and will not present the details here.
\qed

Furthermore, we can easily prove the following proposition, which describes the relations between the variables $T_{1},T_{2}, T_{3}$ and $T_{4}$.
\begin{prop}
The following identities hold:
\begin{enumerate}
\item $T_{1}T_{2}=s^{2}T_{2}T_{1}$;\\
\item $T_{1}T_{3}=r^{2}s^{2}T_{3}T_{1}$;\\
\item $T_{1}T_{4}=r^{2}T_{4}T_{2}$;\\
\item $T_{2}T_{3}=rsT_{3}T_{2}$;\\
\item $T_{2}T_{4}= s^{2}T_{4}T_{2}$;\\
\item $T_{3}T_{4}=rsT_{4}T_{3}$.
\end{enumerate}
\end{prop}
\qed

Now let us denote by $B^{4}$ the subalgebra of $Q(\V)$ generated by $T_{1}^{\pm
1}, X_{2}, X_{3}, X_{4}$, then we have the following
\begin{prop}
The subalgebra $B^{4}$ is the same as the subalgebra of $Q(\V)$ generated
by $X_{1}^{\pm 1}, X_{2}, X_{3}, X_{4}$. In particular, $B^{4}$ is a free module over the subalgebra generated by $X_{2}, X_{3}, X_{4}$.
\end{prop}
\qed

Furthermore, let us denote by $B^{3}$ the subalgebra of $Q(\V)$ generated by $T_{1}^{\pm 1}, T_{2}^{\pm 1}, T_{3}, T_{4}$. Then we shall have the following proposition
\begin{prop}
The subalgebra $B^{3}$ is the same as the subalgebra of $Q(\V)$ generated
by $X_{1}^{\pm 1}, X_{2}^{\pm 1}, X_{3}, X_{4}$.
\end{prop}
\qed

In addition, we will denote by $B^{2}$ the subalgebra of $Q(\U)$ generated by $T_{1}^{\pm 1}, T_{2}^{\pm 1}, T_{3}^{\pm 1}, T_{4}$. We denote by $B^{1}$ the subalgebra of $Q(\U)$ generated by the variables $T_{1}^{\pm 1}, T_{2}^{\pm 1}, T_{3}^{\pm 1}, T_{4}^{\pm 1}$. It is easy to see that the indeterminates $T_{1}, T_{2}, T_{3}, T_{4}$ generate a quantum torus, which we shall denote by 
\[
Q_{4}=\C_{r,s}[T_{1}^{\pm 1}, T_{2}^{\pm 1}, T_{3}^{\pm }, T_{4}^{\pm 1}].
\]

In particular, we have the following proposition:
\begin{prop}
The algebra $Q_{4}=\C_{r,s}[T_{1}^{\pm 1}, T_{2}^{\pm 1},T_{3}^{\pm 1}, T_{4}^{\pm 1}]$ is a quantum torus.
\end{prop}
\qed

Now let us define a linear map 
\[
\mathcal{I} \colon \U \longrightarrow Q_{4}
\]
from $\U$ into $B^{1}=Q_{4}$, whose effect on the generators is given as follows
\begin{eqnarray*}
\mathcal{I}(X_{1})=T_{1}, \quad
\mathcal{I}(X_{2})=T_{2},\quad
\mathcal{I}(X_{3})=T_{3},\\
\mathcal{I}(X_{4})=\lambda (T_{4}+(s^{4}-r^{2}s^{2})T_{2}^{-1}T_{3}+(r^{-1}s-1)T_{2}T_{1}^{-1})
\end{eqnarray*}
where $\lambda=\frac{r}{(r^{2}-s^{2})(r-s)}$.

It is easy to see that the linear map $\mathcal{I}$ can be extended to an algebra monomorphism from $\U$ into the quantum torus $B^{1}=Q_{4}$. Furthermore, it is straightforward to prove the following result:
\begin{thm}
Let us set $B^{5}=\U$ and the following
\begin{eqnarray*}
\Sigma_{5}=\{T_{1}^{i} \mid i \in \Z_{\geq 0}\},\quad
\Sigma_{4}=\{T_{2}^{i}\mid i \in \Z_{\geq 0}\},\\ 
\Sigma_{3}=\{T_{3}^{i}\mid i\in \Z_{\geq 0}\},\quad \Sigma_{2}=\{T_{4}^{i}\mid i\in \Z_{\geq 0}\},
\end{eqnarray*}
then we have the following
\begin{enumerate}
\item $B^{4}=B^{5}\Sigma_{5}^{-1}$;\\
\item $B^{3}=B^{4}\Sigma_{4}^{-1}$;\\
\item $B^{2}=B^{3}\Sigma_{3}^{-1}$;\\
\item $B^{1}=B^{2}\Sigma_{2}^{-1}$;\\
\item The center of $B^{i}$ is the base field $\C$ for $i=1,2,3, 4, 5$.
\end{enumerate}
\end{thm}
\qed

Thanks to the result in \cite{OP}, one knows that any derivation $D$ of
the quantum torus $B^{1}=Q_{4}$ is of the form $D=ad_{t}+\delta$, 
where $ad_{t}$ is an inner derivation defined by some element $t\in B^{1}$, 
and $\delta$ is a central derivation which acts on the variables 
$T_{i}, i=1, 2, 3, 4$ as follows:
\[
\delta(T_{i})=\alpha_{i} T_{i}
\]
for $\alpha_{i}\in \C$.

Suppose that $D$ is a derivation of the algebra $\U=B^{5}$. Due to the nature of the algebras $B^{4}, B^{3}, B^{2}$ and $B^{1}$, we can extend the derivation $D$ to a derivation of the algebras $B^{4}, B^{3}, B^{2}$ and $B^{1}$ respectively.  We will still denote the extended derivations by $D$. Since $B^{1}=Q_{4}$ is a quantum torus and $D$ is derivation of $B^{1}$, we have the following decomposition
\[
D=ad_{t}+\delta
\]
where $ad_{t}$ is an inner derivation determined by some element $t\in B^{1}$,
and $\delta$ is a central derivation of $B^{1}$, which is defined by $\delta(T_{i})=\alpha_{i}T_{i}$ for $\alpha_{i} \in \C, i=1,2,3, 4$. 

We shall prove that the element $t$ can actually be chosen from the algebra $\U=B^{5}$ and the scalars $\alpha_{1},\alpha_{2},\alpha_{3}, \alpha_{4}$ are somehow related to each other. In particular, we shall prove the following key lemma.
\begin{lem}
The following statements are true:
\begin{enumerate}
\item The element $t$ can be chosen from $\U$;\\
\item We have $\alpha_{2}=\alpha_{1}+\alpha_{4}$;\\
\item We have $\alpha_{3}=\alpha_{1}+\alpha_{4}$;\\
\item We have $D(X_{i})=ad_{t}(X_{i})+\alpha_{i}X_{i}$ for $i=1,2, 3, 4$.\\
\end{enumerate}
\end{lem} 
{\bf Proof:} We start the proof by showing that the element $t$ can actually be chosen from the algebra $B^{2}$. Suppose that we have  the following expression of $t$ in the algebra $B^{1}$:
\[
t=\sum_{i,j, k, l}\lambda_{i,j,k, l}T_{1}^{i}T_{2}^{j}T_{3}^{k}T_{4}^{l}. 
\]
If the index $l\geq 0$ for all $l$, then we have proved $t\in B^{2}$. Otherwise, let us set
two elements 
\[
t_{-}=\sum_{l<0}\lambda_{i,j,k, l }T_{1}^{i}T_{2}^{j}T_{3}^{k}T_{4}^{l}
\]
and 
\[
t_{+}=\sum_{l\geq 0}\lambda_{i,j,k, l}T_{1}^{i}T_{2}^{j}T_{3}^{k}T_{4}^{l}.
\]
such that we have $t=t_{-}+t_{+}$. 

Since $D$ is a derivation of the algebra $B^{1}$ and $T_{1} \in B^{1}$, we can apply the derivation $D$ to $T_{1}$. And we obtain the following
\begin{eqnarray*}
D(T_{1}) &=& ad_{t}(T_{1})+\delta(T_{1})\\
&=& (t_{-}T_{1}-T_{1}t_{-})+(t_{+}T_{1}-T_{1}t_{+})+\alpha_{1}(T_{1})
\end{eqnarray*}
for some $\alpha_{1} \in \C$.

Since $D$ is also regarded as a derivation of the algebra $\U$ and the variable $T_{1}$ is also in the algebra $B^{5}=\U$, we shall have that the element $D(T_{1})$ is also in the algebra $B^{5}$, and furthermore in the algebra $B^{2}$. Note that all the elements of $B^{2}$ don't involve negative powers of the variable $T_{4}$, thus we shall have the following
\[
t_{-}T_{1}-T_{1}t_{-}=0. 
\]

Therefore, we are supposed to have the following
\begin{eqnarray*}
T_{1}(\sum_{l<0}\lambda_{i,j,k, l}T_{1}^{i}T_{2}^{j}T_{3}^{k}T_{4}^{l})&=& 
(\sum_{l<0}\lambda_{i,j,k, l, m, n}r^{2k+2l}s^{2j}T_{1}^{i}\\
&&T_{2}^{j}T_{3}^{k}T_{4}^{l})T_{1}\\
&=& (\sum_{k<0}a_{i,j,k, l}T_{1}^{i}T_{2}^{j}T_{3}^{k}T_{4}^{l})T_{1}.
\end{eqnarray*}

This shows that we have the following equations:
\begin{eqnarray*}
 2k+2l=0;\\
2j+2k=0.
\end{eqnarray*}

In addition, applying $D$ to $T_{2}$, we can further derive the following equations
\begin{eqnarray*}
k-2i=0;\\
k+2l=0.
\end{eqnarray*}

Together, these equations show that we have $i=j=k=l=0$, which is a contradiction. Therefore, we have that $t^{-}=0$, which implies $t\in B^{2}$. A similar argument can be used to prove that $t\in B^{3}$.

Since the algebra $B^{3}$ is also generated by the elements $T_{1}^{\pm 1}, T_{2}^{\pm 1},
X_{3}, X_{4}$, we have the following decomposition of $t$ in $B^{3}$: 
\[
t=\sum_{i,j, k\geq 0, l\geq
0}\lambda_{i,j,k,l}T_{1}^{i}T_{2}^{j}X_{3}^{k}X_{4}^{l}.
\]

Applying the derivation $D$ to the variable $T_{1}=X_{1}$, we can further
prove that $j=l$, which implies that $j\geq 0$. Therefore, we 
have proved that $t\in B^{4}$ as desired. Using a similar argument, we can prove that $t\in B^{5}=\U$ as desired.

Since we have $D=ad_{t}+\delta$ for some $t\in \U$ and some central derivation of $B^{1}$, and $T_{1}=X_{1}, T_{2}=X_{2}$ and $T_{3}=X_{3}$, we have  the following
\[
D(X_{i})=D(T_{i})=ad_{t}T_{i}+\alpha_{i}T_{i}=ad_{t}X_{i}+\alpha_{i}X_{i}
\]
for $i=1, 2, 3$.

In addition, we have the following
\begin{eqnarray*}
D(X_{4})&=&ad_{t}(X_{4})+\lambda \delta(T_{4}+(s^{4}-r^{2}s^{2})T_{2}^{-1}T_{3}+(r^{-1}s-1)T_{2}T_{1}^{-1})\\
&=& (tX_{4}-X_{4}t)+\lambda \alpha_{4}T_{4}+\lambda (\alpha_{3}-\alpha_{2})(s^{4}-r^{2}s^{2})T_{2}^{-1}T_{3}\\
&&+\lambda (\alpha_{2}-\alpha_{1})(r^{-1}s-1)T_{2}T_{1}^{-1}\\
&=& (tX_{4}-X_{4}t)+\alpha_{4}X_{4}+\lambda (\alpha_{3}-\alpha_{2}-\alpha_{4})(s^{4}-r^{2}s^{2})T_{2}^{-1}T_{3}\\
&&+\lambda (\alpha_{2}-\alpha_{1}-\alpha_{4})(r^{-1}s-1)T_{2}T_{1}^{-1}.
\end{eqnarray*}

Since the element $D(X_{4})$ is in the algebra $\U$, we shall have
\[
(\alpha_{3}-\alpha_{2}-\alpha_{4})(s^{4}-r^{2}s^{2})T_{2}^{-1}T_{3}+(\alpha_{2}-\alpha_{1}-\alpha_{4})(r^{-1}s-1)T_{2}T_{1}^{-1}=0.
\]
Thus we shall have the following $\alpha_{3}=\alpha_{1}+2\alpha_{4}$ and  $ \alpha_{2}=\alpha_{1}+\alpha_{4}$. In particular, we have
\[
D(X_{4}))=ad_{t}(X_{4})+\alpha_{4}X_{4}.
\]

So we have the proved the lemma as desired.
\qed

Now let us define two derivations $D_{1}, D_{2}$ of the algebra $\U$ 
as follows:
\begin{eqnarray*}
D_{1}(X_{1})=X_{1},\quad D_{X}(X_{2})=X_{2},\quad D_{1}(X_{3})=X_{3},\quad D_{1}(X_{4})=0;\\
D_{2}(X_{1})=0,\quad D_{2}(X_{2})=X_{2},\quad D_{2}(X_{3})=X_{3},\quad D_{2}(X_{4})=X_{4}.
\end{eqnarray*}

Thanks to the previous lemma, we can derive the following result
\begin{thm}
Let $D$ be a derivation of $\U$. Then we have 
\[
D=ad_{t}+\mu_{1}D_{1}+\mu_{2}D_{2}
\]
for some $t\in \U$ and $\mu_{i}\in \C$ for $i=1,2$.
\end{thm}
\qed

Recall that the Hochschild cohomology group in degree 1 of $\U$ is
denoted by $HH^{1}(\U)$, which is defined as follows
\[
HH^{1}(\U)\colon = Der(\U)/InnDer(\U).
\]
where $InnDer(\U)\colon= \{ad_{t} \mid t \in \U \}$ is the Lie algebra of 
inner derivations of $\U$. And it is well known that $HH^{1}(\U)$ is a 
module over $HH^{0}(\U) \colon= Z(\U)=\C$. 

We describe the structural of the first degree Hochschild cohomology
group of $\U$ as a vector space over $\C$. In particular, we have the following theorem
\begin{thm}
The following is true:
\begin{enumerate}
\item Every derivation $D$ of $\U$ can be uniquely written as follows:
\[
D=ad_{t}+\mu_{1}D_{1}+\mu_{2}D_{2}
\]
where $ad_{t}\in InnDer(\U)$ and $\mu_{1}, \mu_{2}\in \C$.\\
\item The first Hochschild cohomology group $HH^{1}(\U)$ of $\U$ is a 
two-dimensional vector space spanned by $\overline{D_{1}}$ 
and $\overline{D_{2}}$.
\end{enumerate}
\end{thm}
{\bf Proof:} Suppose that $ad_{t}+\mu_{1}D_{1}+\mu_{2}D_{2}=0$ as a derivation, to finish the proof, we need to show that $\mu_{1}=\mu_{2}=ad_{t}=0$. Let us set a derivation $\delta=\mu_{1}D_{1}+\mu_{2}D_{2}$. Then $\delta$ can regarded as a derivation of the algebra $\U$, which can be further extended to a derivation of $B^{1}$. As a derivation of $B^{1}$, we also have that $ad_{t}+\delta=0$. Besides, we also have the following
\[
\delta(T_{1})=\mu_{1}T_{1},\quad \delta(T_{2})=\mu_{2}T_{2}, \quad \delta(T_{3})=(\mu_{1}+\mu_{2})T_{3}.
\]

Thus the derivation $\delta$ is verified to a central derivation of the quantum
torus $B^{1}$. Therefore, according to the result in \cite{OP}, we shall have 
that $ad_{t}=0=\delta$. Hence we have $\mu_{1}=\mu_{2}=0$ as desired. 
This proves the uniqueness of the decomposition of the derivation $D$, which 
further proves the second statement of the theorem.
\qed

\section{Hopf algebra automorphisms of the Hopf algebra $\V$}
In this section, we will first determine the algebra automorphism group of the Hopf algebra $\V$. As a result, we are able to determine the Hopf algebra automorphism group of $\V$ as well. We will closely follow the approach used in \cite{F}. Note that such an approach has been adopted to investigate the automorphism group of $\check{U}_{r,s}^{\geq 0}(\mathfrak sl_{3})$ in \cite{T1}. It is no surprise that we derive very similar results to those obtained in \cite{T1}.
\subsection{The Hopf algebra $U_{r,s}^{\geq 0}(B_{2})$}
To proceed, we first recall the construction of the  Hopf subalgebra of $U_{r,s}^{\geq 0}(B_{2})$ from \cite{BGH}. Later on, we will define an augmented version of the Hopf algebra $U_{r,s}^{\geq 0}(B_{2})$. 
\begin{defn}
The Hopf algebra $U_{r,s}^{\geq 0}(B_{2})$ is defined to be the $\C-$algebra generated by the generators $e_{1}, e_{2}$ and $w_{1}, w_{2}$ subject to the following relations:
\begin{eqnarray*}
w_{1}w_{1}^{-1}=w_{2}w_{2}^{-1}=1,\quad w_{1}w_{2}=w_{2}w_{1};\\
w_{1}e_{1}=r^{2}s^{-2}e_{1}w_{1},\quad w_{1}e_{2}=s^{2}e_{2}w_{1};\\
w_{2}e_{1}=r^{-2}e_{1}w_{2},\quad w_{2}e_{2}=rs^{-1}e_{2}w_{2};\\
e_{1}^{2}e_{2}-(r^{2}+s^{2})e_{1}e_{2}e_{1}+r^{2}s^{2}e_{2}e_{1}^{2}=0;\\
e_{1}e_{2}^{3}-(r^{2}+rs+s^{2})e_{2}e_{1}e_{2}^{2}
+rs(r^{2}+rs+s^{2})e_{2}^{2}e_{1}e_{2}-r^{3}s^{3}e_{2}^{3}e_{1}=0.
\end{eqnarray*}
\end{defn}

It can be easily verified that the following operators define a Hopf algebra structure on $U^{\geq 0}_{r,s}(B_{2})$.
\begin{eqnarray*}
\Delta(e_{1})=e_{1}\otimes 1+w_{1}\otimes e_{1};\\
\Delta(e_{2})=e_{2}\otimes 1+w_{2}\otimes e_{2};\\
\Delta(w_{1})=w_{1}\otimes w_{1},\quad \Delta(w_{2})=w_{2}\otimes w_{2};\\
S(e_{1})=-w_{1}e_{1},\quad S(e_{2})=-w_{2}e_{2};\\
S(w_{1})=w_{1}^{-1},\quad S(w_{2})=w_{2}^{-1};\\
\epsilon(e_{1})=\epsilon(e_{2})=0,\quad \epsilon (w_{1})=\epsilon(w_{2})=1.
\end{eqnarray*}

Of course, it is easy to see that we have the following proposition:
\begin{prop}
The set 
\[
\{X_{1}^{a}X_{2}^{b}X_{3}^{c}X_{4}^{d}w_{1}^{m}w_{2}^{n}|a, b, c, d \in \Z_{\geq 0}, m, n \in \Z\}
\]

forms a PBW-basis of the algebra $U_{r,s}^{\geq 0}(B_{2})$.
\end{prop}
\qed

However, we will not study the algebra $U^{\geq 0}_{r,s}(B_{2})$ in this paper. Instead, in the rest of this paper, we will study  the automorphisms of an augmented Hopf algebra $\V$, whose definition will be given in the next subsection.

\subsection{The Augmented Hopf algebra $\V$}
In this subsection, we shall introduce an augmented Hopf algebra $\V$, which contains the algebra $U^{+}_{r,s}(B_{2})$ as a subalgebra and enlarges the Hopf algebra $U^{\geq 0}_{r,s}(B_{2})$. To this end, we need to define the following new variables:
\[
k_{1}=w_{1}w_{2},\quad k_{2}=w_{1}^{1/2}w_{2}.\\
\]

It is easy to see that 
\[
w_{1}=k_{1}^{2}k_{2}^{-2},\quad w_{2}=k_{1}^{-1}k_{2}^{2}.
\]

Now we can have the following definition of the algebra $\V$.
\begin{defn}
The algebra ${\check U}^{\geq 0}_{r,s}({\mathfrak sl_{3}})$ is a $\C$-algebra generated by 
$e_{1}, e_{2}, k_{1}^{\pm 1},$ and $ k_{2}^{\pm 1}$ subject to the following relations:
\begin{eqnarray*}
k_{1}k_{1}^{-1}=1=k_{2}k_{2}^{-1},\quad k_{1}k_{2}=k_{2}k_{1};\\
k_{1}e_{1}=s^{-2}e_{1}k_{1},\quad k_{1}e_{2}=rs e_{2}k_{1};\\
k_{2}e_{1}=r^{-1}s^{-1}e_{1}k_{2},\quad k_{2}e_{2}=re_{2}k_{2};\\
e_{1}^{2}e_{2}-(r^{2}+s^{2})e_{1}e_{2}e_{1}+r^{2}s^{2}e_{2}e_{1}^{2}=0;\\
e_{1}e_{2}^{3}-(r^{2}+rs+s^{2})e_{2}e_{1}e_{2}^{2}
+rs(r^{2}+rs+s^{2})e_{2}^{2}e_{1}e_{2}-r^{3}s^{3}e_{2}^{3}e_{1}=0.
\end{eqnarray*}
\end{defn}
\qed

We now introduce a Hopf algebra structure on $\V$ by defining the following operators:
\begin{eqnarray*}
\Delta(e_{1})=e_{1}\otimes 1+k_{1}^{2}K_{2}^{-2}\otimes e_{1};\\
\Delta(e_{2})=e_{2}\otimes 1+k_{1}^{-1}K_{2}^{2}\otimes e_{2};\\
\Delta(k_{1})=k_{1}\otimes k_{1},\quad \Delta(k_{2})=k_{2}\otimes k_{2};\\
S(e_{1})=-k_{1}^{2}k_{2}^{-2}e_{1},\quad S(e_{2})=-K_{1}^{-1}k_{2}^{2}e_{2};\\
S(k_{1})=k_{1}^{-1},\quad 
S(k_{2})=k_{2}^{-1};\\
\epsilon(e_{1})=\epsilon(e_{2})=0,\quad
\epsilon (k_{1})=\epsilon(k_{2})=1.
\end{eqnarray*}

Then it is straightforward to verify the following result:
\begin{prop}
The algebra ${\check U}^{\geq 0}_{r,s}({\mathfrak sl_{3}})$ is a Hopf
algebra with the coproduct, counit and antipode defined as above.
\end{prop}
\qed

Furthermore, it is easy to see that we have the following result 
\begin{thm}
The Hopf algebra ${\check U}^{\geq 0}_{r,s}(B_{2})$ has a $\C-$basis 
\[
\{ k_{1}^{m}k_{2}^{n}X_{1}^{a}X_{2}^{b}X_{3}^{c}X_{4}^{d}\mid m, n\in \Z,\, a,b,c,d \in \Z_{\geq 0}\}.
\]
\end{thm}
\qed

In particular, one can see that all the invertible elements of $\V$ are of the form $\lambda k_{1}^{m}k_{2}^{n}$ for some $\lambda \in \C^{\ast}$ and $m, n \in \Z$.
\subsection{The algebra automorphism group of $\V$}

Suppose that $\theta$ denotes an algebra automorphism of the Hopf algebra $\V$. Since $k_{1}, k_{2}$ are invertible elements in the algebra $\V$ and $\theta$ is 
an algebra automorphism of $\V$, the images $\theta(k_{1}), \theta(k_{2})$  of $k_{1}, k_{2}$ are invertible elements in $\V$. Recall that the invertible elements of 
the algebra $\V$ are of the form $\lambda k_{1}^{m}k_{2}^{n}, \lambda \in \C^{\ast}, m, n \in \Z$. Therefore, we shall have the 
following expressions of $\theta(k_{1})$ and $\theta(k_{2})$: 
\[
\theta(k_{1})=\lambda_{1}k_{1}^{x}k_{2}^{y},\quad \theta(k_{2})=\lambda_{2}k_{1}^{z}k_{2}^{w}
\]
where $\lambda_{1}, \lambda_{2}\in \C^{\ast}$ and $x,y,z,w \in \Z$. 

Since $\theta$ is an algebra automorphism of $\V$, there is an invertible $2\times 2$ matrix associated to $\theta$. We will denote this matrix by $M_{\theta}=(M_{ij})$. 
As a matter of fact, we will set the entries $M_{11}=x, M_{12}=y, M_{21}=z$ and $M_{22}=w$. Since $\theta$ is an algebra automorphism, we know the determinant 
of $M_{\theta}$ is $\pm 1$, i.e., we have $xw-yz=\pm 1$.

For $l=1,2$, let us set the following expressions of the images of $e_{1}, e_{2}$ under the automorphism $\theta$:
\[
\theta(e_{l})=\sum_{m_{l}, n_{l}, \beta_{l}^{1}, \beta_{l}^{2}, \beta_{l}^{3},\beta_{l}^{4}}\gamma_{m_{l}, n_{l}, \beta_{l}^{1}, \beta_{l}^{2}, \beta_{l}^{3}, \beta_{l}^{4}}k_{1}^{m_{l}}k_{2}^{n_{l}}X_{1}^{\beta_{l}^{1}}X_{2}^{\beta_{l}^{2}}X_{3}^{\beta_{l}^{3}} X_{4}^{\beta_{l}^{4}}
\]
where $\gamma_{m_{l}, n_{l},\beta_{l}^{1}, \beta_{l}^{2}, \beta_{l}^{3}, \beta_{l}^{4}} \in \C^{\ast}$ and
$m_{l}, n_{l} \in \Z$ and $\beta_{l}^{1}, \beta_{l}^{2}, \beta_{l}^{3}, \beta_{l}^{4} \in \Z_{\geq 0}$. 

Then we have the following
\begin{prop}
Suppose that $\theta \in Aut_{\C}(\V)$ is an algebra automorphism of $\V$, then we have $M_{\theta} \in GL(2,\Z_{\geq 0})$.
\end{prop}
{\bf Proof:} Since $k_{1}e_{1}=s^{-2}e_{1}k_{1}$ and $k_{2}e_{1}=r^{-1}s^{-1}e_{1}k_{2}$ and $\theta$ is an algebra automorphism of $\V$, we have the following
\begin{eqnarray*}
\theta(k_{1})\theta(e_{1})=s^{-2}\theta(e_{1})\theta(k_{1});\\
\theta(k_{2})\theta(e_{1})=r^{-1}s^{-1}\theta(e_{1})\theta(k_{2});
\end{eqnarray*}
which implies that
\begin{eqnarray*}
& &\lambda_{1}k_{1}^{x}k_{2}^{y}(\sum_{m_{1}, n_{1}, \beta_{1}^{1}, \beta_{1}^{2}, \beta_{1}^{3},\beta_{1}^{4}}\gamma_{m_{1}, n_{1}, \beta_{1}^{1}, \beta_{1}^{2}, \beta_{1}^{3}, \beta_{1}^{4}}k_{1}^{m_{1}}k_{2}^{n_{1}}X_{1}^{\beta_{1}^{1}}X_{2}^{\beta_{1}^{2}}X_{3}^{\beta_{l}^{3}} X_{4}^{\beta_{1}^{4}})\\
&=& s^{-2} (\sum_{m_{1}, n_{1}, \beta_{1}^{1}, \beta_{1}^{2}, \beta_{1}^{3},\beta_{1}^{4}}\gamma_{m_{1}, n_{1}, \beta_{1}^{1}, \beta_{1}^{2}, \beta_{1}^{3}, \beta_{1}^{4}}k_{1}^{m_{1}}k_{2}^{n_{1}}X_{1}^{\beta_{1}^{1}}X_{2}^{\beta_{1}^{2}}X_{3}^{\beta_{1}^{3}} X_{4}^{\beta_{1}^{4}})\lambda_{1}k_{1}^{x}k_{2}^{y},
\end{eqnarray*}
and 
\begin{eqnarray*}
& & \lambda_{2}k_{1}^{z}k_{2}^{w}(\sum_{m_{1}, n_{1}, \beta_{1}^{1}, \beta_{1}^{2}, \beta_{1}^{3},\beta_{1}^{4}}\gamma_{m_{1}, n_{1}, \beta_{1}^{1}, \beta_{1}^{2}, \beta_{1}^{3}, \beta_{1}^{4}}k_{1}^{m_{l}}k_{2}^{n_{l}}X_{1}^{\beta_{l}^{1}}X_{2}^{\beta_{l}^{2}}X_{3}^{\beta_{l}^{3}} X_{4}^{\beta_{l}^{4}}) \\
&=& r^{-1}s^{-1}(\sum_{m_{l}, n_{l}, \beta_{l}^{1}, \beta_{l}^{2}, \beta_{l}^{3},\beta_{1}^{4}}\gamma_{m_{1}, n_{1}, \beta_{1}^{1}, \beta_{l}^{2}, \beta_{l}^{3}, \beta_{1}^{4}}k_{1}^{m_{1}}k_{2}^{n_{1}}X_{1}^{\beta_{1}^{1}}X_{2}^{\beta_{1}^{2}}X_{3}^{\beta_{l}^{3}} X_{4}^{\beta_{1}^{4}})\lambda_{2}k_{1}^{z}k_{2}^{w}.
\end{eqnarray*}

After the calculations,  we shall have the following system of equations:
\begin{eqnarray*}
(2\beta_{1}^{1}+2\beta_{1}^{2}+2\beta_{1}^{3})x+(\beta_{1}^{2}+2\beta_{1}^{3}+\beta_{1}^{4})y=2;\\
(2\beta_{1}^{1}+2\beta_{1}^{2}+2\beta_{1}^{3})z+(\beta_{1}^{2}+2\beta_{1}^{3}+\beta_{1}^{4})w=0.
\end{eqnarray*}

Similarly, we also have the following
\begin{eqnarray*}
(\beta_{2}^{1}+\beta_{2}^{2}+\beta_{2}^{3})x+(\beta_{2}^{2}+2\beta_{2}^{3}+\beta_{2}^{4})y=0;\\
(\beta_{2}^{1}+\beta_{2}^{2}+\beta_{2}^{3})z+(\beta_{2}^{2}+2\beta_{2}^{3}+\beta_{2}^{4})w=1.
\end{eqnarray*}

We now define a $2\times2-$matrix $B=(b_{ij})$ with the following entries
\begin{eqnarray*}
b_{11}=2\beta_{1}^{1}+2\beta_{1}^{2}+2\beta_{1}^{3};\\
b_{21}=\beta_{1}^{2}+2\beta_{1}^{3}+\beta_{1}^{4};\\
b_{12}=\beta_{2}^{1}+\beta_{2}^{2}+\beta_{2}^{3};\\
b_{22}=\beta_{2}^{2}+2\beta_{2}^{3}+\beta_{2}^{4}.
\end{eqnarray*}

And we shall have the following system of equations
\begin{eqnarray*}
b_{11}x+b_{21}y=2;\\
b_{12}x+b_{22}y=0;\\
b_{11}z+b_{21}w=0;\\
b_{12}z+b_{22}w=1.
\end{eqnarray*}

This system of equations implies that we have the following
\[
M_{\theta}B=\left( \begin{array}{lr}
2 & 0\\
0 &1 
\end{array}
\right)
\]
which shows that we have 
\[
M_{\theta}^{-1}=\left( \begin{array}{lr}
b_{11}/2 & b_{12}\\
b_{21}/2 &b_{22} 
\end{array}
\right)
\] 

Obviously, we have that $M_{\theta}^{-1}=M_{\theta^{-1}}$, where the matrix $M_{\theta^{-1}}$ is the corresponding matrix associated to 
the algebra automorphism $\theta^{-1}$. Since the entries $b_{11}, b_{12}, b_{21}, b_{22}$ are all nonnegative integers, we can conclude 
that $M_{\theta^{-1}}\in GL(2, \Z_{\geq 0})$. Applying the similar arguments to the algebra automorphism $\theta^{-1}$, we shall be able 
prove that $M_{\theta}\in GL(2, \Z_{\geq 0})$ as desired. 
\qed

For the reader's convenience, we recall an important lemma ({\bf Lemma 2.2} from \cite{F}), which characterizes 
the matrix $M_{\theta}$:
\begin{lem}
If $M$ is a matrix in $GL(n,\Z_{\geq 0})$ such that it’s inverse matrix $M^{-1}$ is also in $GL(n,\Z_{\geq 0})$, then we 
have $M=(\delta_{i\sigma(j)})_{i,j}$, where $\sigma$ is an element of the symmetric group $\mathbb{S}_{n}$.
\end{lem}
\qed

As a result of  {\bf Proposition 2.3} and {\bf Lemma 2.1}, we immediately have the following result
which describes the images of $k_{1}, k_{2}$ under an automorphism $\theta$ of $\V$:
\begin{cor}
Let $\theta \in Aut_{\C}(\V)$ be an algebra automorphism of $\V$. Then for $l=1,2$, we have 
\[
\theta(k_{l})=\lambda_{l} k_{\sigma(l)}
\]
where $\sigma \in \mathbb{S}_{2}$ and $\lambda_{l} \in \C^{\ast}$.
\end{cor}
\qed

Furthermore, we can prove the following result:
\begin{prop}
Let $\theta \in Aut_{\C}(\V)$ be an algebra automorphism of $\V$. Then for $l=1,2$, we have 
\[
\theta(e_{l})=\gamma_{l}k_{1}^{m_{l}}k_{2}^{n_{l}}e_{\sigma(l)}
\]
where $\gamma_{l}\in \C^{\ast}$ and $m_{l}, n_{l} \in \Z$.
\end{prop}
{\bf Proof:} Let $\theta \in Aut_{\C}(\V)$ be an algebra automorphism
of $\V$. We will need to consider two cases.

{\bf Case 1}: Suppose that $\theta(k_{1})=\lambda_{1}k_{1}$ and $\theta(k_{2})=\lambda_{2}k_{2}$, then it suffices 
to show that we have  
\[
\theta(e_{1})=\gamma_{1}k_{1}^{m_{1}}k_{2}^{n_{1}}e_{1},
\quad
\theta(e_{2})=\gamma_{2}k_{1}^{m_{2}}k_{2}^{n_{2}}e_{2}
\]
for some $\gamma_{1}, \gamma_{2} \in \C^{\ast}$ and $m_{1},m_{2}, n_{1}, n_{2} \in \Z$.

Note that we have the following relations between $e_{1}, e_{2}$ and $k_{1}, k_{2}$:
\[
k_{1}e_{1}=s^{-2}e_{1}k_{1},\quad k_{2}e_{1}=r^{-1}s^{-1}e_{1}k_{2}.
\]

Via applying $\theta$ to these identities, we shall have the  following
\begin{eqnarray*}
\theta(k_{1})\theta(e_{1})=s^{-2}\theta(e_{1})\theta(k_{1});\\
\theta(K_{2})\theta(e_{1})=r^{-1}s^{-1}\theta(e_{2})\theta(k_{1}).
\end{eqnarray*}

Therefore, we shall the following
\begin{eqnarray*}
&&\lambda_{1}k_{1}(\sum_{m_{1},n_{1},\beta_{1}^{1},\beta_{1}^{2},
\beta_{1}^{3}, \beta_{1}^{4}}\gamma_{m_{1}, n_{1}, \beta_{1}^{1}, \beta_{1}^{2}, \beta_{1}^{3},\beta_{1}^{4}}
k_{1}^{m_{1}}k_{2}^{n_{1}}X_{1}^{\beta_{1}^{1}}X_{2}^{\beta_{1}^{2}}X_{3}^{\beta_{1}^{3}}X_{4}^{\beta_{1}^{4}})\\
&=&\lambda_{1}s^{-2}(\sum_{m_{1},
n_{1}, \beta_{1}^{1},\beta_{1}^{2},\beta_{1}^{3}}\gamma_{m_{1}, n_{1},
\beta_{1}^{1},\beta_{1}^{2},
\beta_{1}^{3}}K_{1}^{m_{1}}k_{2}^{n_{1}}X_{1}^{\beta_{1}^{1}}X_{2}^{\beta_{1}^{2}}X_{3}^{\beta_{1}^{3}}
X_{4}^{\beta_{1}^{4}})k_{1}.
\end{eqnarray*}

In addition, we also have the following
\begin{eqnarray*}
&&\lambda_{2}k_{2}(\sum_{m_{1},n_{1},\beta_{1}^{1},\beta_{1}^{2},
\beta_{1}^{3}}\gamma_{m_{1}, n_{1}, \beta_{1}^{1}, \beta_{1}^{2}, \beta_{2}^{3},\beta_{2}^{4}}
k_{1}^{m_{1}}k_{2}^{n_{1}}X_{1}^{\beta_{1}^{1}}X_{2}^{\beta_{1}^{2}}X_{3}^{\beta_{1}^{3}}X_{4}^{\beta_{1}^{4}})\\
&=&\lambda_{2}r^{-1}s^{-1}(\sum_{m_{1},n_{1}, \beta_{1}^{1},\beta_{1}^{2},\beta_{1}^{3}}\gamma_{m_{1}, n_{1},
\beta_{1}^{1},\beta_{1}^{2},
\beta_{1}^{3},\beta_{1}^{4}}k_{1}^{m_{1}}K_{2}^{n_{1}}X_{1}^{\beta_{1}^{1}}
X_{2}^{\beta_{1}^{2}}X_{3}^{\beta_{1}^{3}}X_{4}^{\beta_{1}^{4}})k_{2}.
\end{eqnarray*}

Thus we shall have the following
\begin{eqnarray*}
s^{-(\beta_{1}^{1}+\beta_{1}^{2}+\beta_{1}^{3})}(rs)^{(\beta_{1}^{2}+2\beta_{1}^{3}+\beta_{1}^{4})}=s^{-2},\\
(rs)^{-(\beta_{1}^{1}+\beta_{1}^{2}+\beta_{1}^{3})}
r^{(\beta_{1}^{2}+2\beta_{1}^{3}+\beta_{1}^{4})}=r^{-1}s^{-1}.
\end{eqnarray*}

Moreover, the above identities imply the following
\begin{eqnarray*}
2\beta_{1}^{1}+2\beta_{1}^{2}+2\beta_{1}^{3}=2;\\
\beta_{1}^{2}+2\beta_{1}^{3}+\beta_{1}^{4}=0;\\
\beta_{1}^{1}+\beta_{1}^{2}+\beta_{1}^{3}=1;\\
\beta_{1}^{2}+2\beta_{1}^{3}+\beta_{1}^{4}=0.
\end{eqnarray*}

Note that all $\beta_{i}^{j}, i, j =1, 2, 3, 4$ are non-negative integers, thus we shall have that 
\[
\beta_{1}^{1}=1,\quad \beta_{1}^{2}=\beta_{1}^{3}=\beta_{1}^{4}=0.
\]

Similarly, we can also verify the following
\[
\beta_{2}^{1}=\beta_{2}^{2}=\beta_{2}^{3}=0,\quad \beta_{2}^{4}=1.
\] 
Therefore, we have proved the result for {\bf Case 1}.

{\bf Case 2}: Suppose that $\theta(k_{1})=\lambda_{1}k_{2}$ and $\theta(k_{2})=\lambda_{2}k_{1}$, we have to 
prove that $\theta(e_{1})=\gamma_{1}k_{1}^{m_{1}}k_{2}^{n_{1}}e_{2}$ and $\theta(e_{2})=\gamma_{2}k_{1}^{m_{2}}k_{2}^{n_{2}}e_{1}$. 
We will not repeat the details here because the proof goes the same way as in {\bf Case 1}.
\qed 

Now we are going to verify that, in a sense, the generators $e_{1}, e_{2}$ can not be exchanged by any algebra automorphism $\theta$ of $\V$. 
Indeed, we have the following result
\begin{cor}
Let $\theta \in Aut_{\C}(\V)$ be an algebra automorphism of $\V$. Then for $l=1,2$, we have the following 
\[
\theta(k_{l})=\lambda_{l}k_{l},\, \theta(e_{l})=\gamma_{l}k_{1}^{m_{l}}k_{2}^{n_{l}}e_{l}
\]
where $\lambda_{l}, \gamma_{l} \in \C^{\ast}$ and $m_{l},n_{l} \in \Z$.
\end{cor}
{\bf Proof:} Suppose that $\theta(k_{1})=\lambda_{1}k_{2}$ and $\theta(e_{2})=\gamma_{1}k_{1}^{m_{1}}k_{2}^{n_{1}}e_{2}$. Since we have $\theta(k_{1})\theta(e_{1})=s^{-2}\theta(e_{1})\theta(k_{1})$, we have the following
\[
\lambda_{1}k_{2}\gamma_{1} k_{1}^{m_{1}}k_{2}^{m_{2}}e_{2}=s^{-2}\gamma_{1} k_{1}^{m_{1}}k_{2}^{m_{2}}e_{2}\lambda_{1}k_{2}.
\]
Note that $k_{2}e_{2}=re_{2}k_{2}$, then we got a contradiction. Therefore, we have proved the statement as desired.
\qed

The following main theorem describes the algebra automorphism group of the algebra $\V$:
\begin{thm}
Let $\theta \in Aut_{\C}(\V)$ be an algebra automorphism of the algebra $\V$. Then for $l=1,2$, we have the following
\[
\theta(k_{l})=\lambda_{l}k_{l},\quad \theta(e_{1})=\gamma_{1}k_{1}^{a}K_{2}^{b}e_{1}, \quad \theta(e_{2})=\gamma_{2}k_{1}^{c}k_{2}^{d}e_{2}
\]
where $\lambda_{l}, \gamma_{l} \in \C^{\ast}$ and $a, b, c, d \in \Z$ such that $b=2c, a+2c+d=0$.
\end{thm}
{\bf Proof:} Let $\theta$ be an algebra automorphism of $\V$ and suppose that 
\[
\theta(e_{1})=\gamma_{1}k_{1}^{a}k_{2}^{b}e_{1},\quad
\theta(e_{2})=\gamma_{2}k_{1}^{c}k_{2}^{d}e_{2}.
\]

Note that we have the following
\begin{eqnarray*}
(k_{1}^{a}k_{2}^{b}e_{1})(k_{1}^{a}k_{2}^{b}e_{1})(k_{1}^{c}k_{2}^{d}e_{2})
&=& (s^{2})^{a} (rs)^{b}(s^{2})^{2c}(rs)^{2d}\\
& & k_{1}^{2a+c}k_{2}^{2b+d}e_{1}^{2}e_{2}\\
&=& r^{(b+2d)}s^{(2a+b+4c+2d)}k_{1}^{2a+c}k_{2}^{2b+d}e_{1}^{2}e_{2};
\end{eqnarray*}

and
\begin{eqnarray*}
(k_{1}^{a}k_{2}^{b}e_{1})(k_{1}^{c}k_{2}^{d}e_{2})(k_{1}^{a}k_{2}^{b}e_{1})
&=& (s^{2})^{c} (rs)^{d}((rs)^{-1})^{a}(s^{2})^{a}(r^{-1})^{b}(rs)^{b}\\
&&k_{1}^{2a+c}k_{2}^{2b+d}e_{1}e_{2}e_{1}\\
&=& r^{(d-a)}s^{(a+b+2c+d)} k_{1}^{2a+c}k_{2}^{2b+d}e_{1}e_{2}e_{1};\\
\end{eqnarray*}

and
\begin{eqnarray*}
(k_{1}^{c}k_{2}^{d}e_{2})
(k_{1}^{a}k_{2}^{b}e_{1})(k_{1}^{a}k_{2}^{b}e_{1})&=&
(r^{-1}s^{-1})^{a} (r^{-1})^{b}(s^{2})^{a}(r^{-1}s^{-1})^{a}(rs)^{b}\\
&&k_{1}^{2a+c}k_{2}^{2b+d}e_{2}e_{1}^{2}\\
&=& r^{(-2a-b)}s^{b}k_{1}^{2a+c}k_{2}^{2b+d}e_{2}e_{1}^{2}.
\end{eqnarray*}

Via applying the automorphism $\theta$ to the first two-parameter quantum Serre relation
\[
e_{1}^{2}e_{2}-(r^{2}+rs+s^{2})e_{1}e_{2}e_{1}+(rs)^{2}e_{2}e_{1}^{2}=0
\]
we shall have the following system of equations
\begin{eqnarray*}
b+2d&=& -a+d;\\
-2a-b&=&-a+d;\\
2a+b+4c+2d&=& a+b+2c+d;\\
b&=&a+b+2c+d.
\end{eqnarray*}

It is easy to see that the previous system of equations is reduced to the following system of equations
\begin{eqnarray*}
a+b+d&=&0;\\
a+2c+d &=&0.
\end{eqnarray*}

In addition, direct calculations yield  the following
\begin{eqnarray*}
 (k_{1}^{c}k_{2}^{d}e_{2})(k_{1}^{c}k_{2}^{d}e_{2})(k_{1}^{c}k_{2}^{d}e_{2})(k_{1}^{a}k_{b}e_{1})&=& r^{(-3a-3b-3c-3d)}s^{(-3a-3c)}\\
&& k_{1}^{3c+a}k_{2}^{3d+b}e_{2}^{3}e_{1};
\end{eqnarray*}

and
\begin{eqnarray*}
(k_{1}^{c}k_{2}^{d}e_{e})(k_{1}^{c}k_{2}^{d}e_{2})(k_{1}^{a}k_{2}^{b}e_{1})(k_{1}^{c}k_{2}^{2}e_{2})&=& r^{(-2a-2b-3c-3d)}s^{(-2a-c+d)}\\
&& k_{1}^{3c+a}k_{2}^{3d+b}e_{2}^{2}e_{1}e_{2};
\end{eqnarray*}

and
\begin{eqnarray*}
(k_{1}^{c}k_{2}^{d}e_{2})(k_{1}^{a}k_{2}^{b}e_{1})(k_{1}^{c}k_{2}^{d}e_{2})(k_{1}^{c}k_{2}^{d}e_{2})
&=& r^{(-a-b-3c-d)}s^{-a+c+2d}\\
&&k_{1}^{3c+a}k_{2}^{3d+b}e_{2}e_{1}e_{2}^{2};
\end{eqnarray*}

and 
\begin{eqnarray*}
 (k_{1}^{a}k_{b}e_{1})(k_{1}^{c}k_{2}^{d}e_{2})(k_{1}^{c}k_{2}^{d}e_{2})(k_{1}^{c}k_{2}^{d}e_{2})&=& r^{(-3c)}s^{(3c+3d)}\\
&& k_{1}^{3c+a}k_{2}^{3d+b}e_{1}e_{2}^{3}.
\end{eqnarray*}

Via applying the automorphism $\theta$ to the second two-parameter quantum Serre relation
\begin{eqnarray*}
e_{2}^{3}e_{1}-(r^{-2}+r^{-1}s^{-1}+s^{-2})e_{2}^{2}e_{1}e_{2}+r^{-1}s^{-1}(r^{-2}+r^{-1}s^{-1}\\
+s^{-2}) e_{2}e_{1}e_{2}^{2}-(rs)^{-3}e_{1}e_{2}^{3}=0
\end{eqnarray*}
we shall have the following system of equations
\begin{eqnarray*}
-3a-3b-3c-3d&=& -2a-2b-3c-3d;\\
-a-b-3c-d&=& -2a-2b-3c-3d;\\
-3c & = &-2a-2b-3c-3d;\\
-3a-3c&=& -2a-c+d;\\
-a+c+2d&=& -2a-c+d;\\
3c+3d&=& -2a-c+d.
\end{eqnarray*}

Therefore, we also have the same system of equations as follows
\begin{eqnarray*}
a+b+d&=&0;\\
a+2c+d &=&0.
\end{eqnarray*}

Solving the system
\begin{eqnarray*}
a+b+d&=&0;\\
a+2c+d &=&0;
\end{eqnarray*}
we have that $b=2c$ and $a+2c+d=0$. Thus we have proved the theorem as desired.
\qed

\subsection{Hopf algebra automorphisms of $\V$}

In this subsection, we further determine all the Hopf algebra automorphisms of the 
Hopf algebra $\V$. Let us denote by $Aut_{Hopf}(\V)$ the group of all Hopf algebra automorphisms of $\V$. 

First of all, we have the following result 
\begin{thm}
Let $\theta \in Aut_{Hopf}(\V)$. Then for $l=1,2$, we have the following
\[
\theta(k_{l})=k_{l},\quad \theta(e_{l})=\gamma_{l}e_{l},
\]
for some $\gamma_{l}\in \C^{\ast}$. In particular, we have 
\[
Aut_{Hopf}(\V)\cong (\C^{\ast})^{2}.
\]
\end{thm}
{\bf Proof:} First of all, let $\theta \in Aut_{Hopf}(\V)$ denote a Hopf algebra
automorphism of $\V$, then we have $\theta \in Aut_{\C}(\V)$. Therefore, 
we shall have the following
\begin{eqnarray*}
\theta(k_{l})=\lambda_{l}k_{l};\\
\theta(e_{1})=\gamma_{1}k_{1}^{a}k_{2}^{b}e_{1};\\
\theta(E_{2})=\gamma_{2}k_{1}^{c}k_{2}^{d}e_{2};
\end{eqnarray*} 
for some $\lambda_{l}, \gamma_{l}\in \C^{\ast}$ for $l=1,2$, and $a,b,c,d \in \Z$ such that $b=2c, 
a+2c+d=0$. 

We want to prove that $\lambda_{l}=1$ for $l=1,2$. Since $\theta$ is a
Hopf algebra automorphism, we shall have the following
\[
(\theta \otimes \theta)(\Delta(k_{l}))=\Delta(\theta(k_{l}))
\]
for $l=1,2$, which imply the following
\[
\lambda_{l}^{2}=\lambda_{l}
\]
for $l=1,2$. Therefore, we have $\lambda_{l}=1$ for $l=1,2$. 

Now we need to prove that $a=b=c=d=0$. First of all, note that we have the following
\begin{eqnarray*}
\Delta(\theta(e_{1}))&=&\Delta(\gamma_{1}k_{1}^{a}k_{2}^{b}e_{1})\\
&=&\Delta(\gamma_{1}k_{1}^{a}k_{2}^{b})\Delta(e_{1})\\
&=&\gamma_{1}(k_{1}^{a}k_{2}^{b}\otimes k_{1}^{a}k_{2}^{b}) (e_{1}\otimes 1+k_{1}^{2}k_{2}^{-2}\otimes e_{1})\\
&=& \gamma_{1}k_{1}^{a}k_{2}^{b}e_{1}\otimes k_{1}^{a}k_{2}^{b}+\gamma_{1}k_{1}^{a}k_{2}^{b}k_{1}^{2}k_{2}^{-2}\otimes k_{1}^{a}k_{2}^{b}e_{1}\\
&=&\theta(e_{1})\otimes k_{1}^{a}k_{2}^{b}+k_{1}^{a}k_{2}^{b}k_{1}^{2}k_{2}^{-2}\otimes \theta(e_{1}).
\end{eqnarray*}

Second of all, note that we also have the following
\begin{eqnarray*}
(\theta\otimes \theta)(\Delta(e_{1}))&=&(\theta\otimes \theta)(e_{1}\otimes 1+k_{1}^{2}k_{2}^{-2}\otimes e_{1})\\
&=& \theta(e_{1})\otimes 1+\theta(k_{1}^{2}k_{2}^{-2})\otimes \theta(e_{1})\\
&=& \theta(e_{1})\otimes 1+ k_{1}^{2}k_{2}^{-2}\otimes \theta(e_{1}).
\end{eqnarray*} 

Since $\Delta(\theta(e_{1}))=(\theta\otimes \theta)\Delta(e_{1})$, we have $a=b=0$. Since $b=2c$ and $a+2c+d=0$, we have $a=b=c=d=0$.

In addition, it is obvious that the algebra automorphism $\theta$ defined by $\theta(k_{l})=k_{l}$ and $\theta(e_{l})=\gamma_{l}e_{l}$ for $l=1, 2$ is 
a Hopf algebra automorphism of $\V$. Thus,  we have proved the theorem.
\qed

{\bf Acknowledgement:} The author would like to thank L. Ben Yakoub for providing a copy of his paper during the preparation of this work.

\end{document}